\documentclass[reqno]{amsart}
\usepackage{graphicx}
\usepackage{amsthm,amssymb,amsfonts}

\usepackage{verbatim}

\usepackage{mathUtil}

\usepackage{xcolor}
\usepackage[outline]{contour}

\usepackage[backend=biber,style=alphabetic,maxnames=50]{biblatex}
\usepackage{multirow}

\usepackage{hyperref}
\hypersetup{
    colorlinks=true,
    linkcolor=blue,
    filecolor=magenta,
    urlcolor=cyan,
    pdftitle={Maximality Principles in Modal Logic and the Axiom of Choice},
    pdfpagemode=FullScreen,
}
\usepackage{todonotes}
\usepackage[capitalize]{cleveref}

\begin{document}

\title[Maximality Principles in Modal Logic and the Axiom of Choice]{Maximality Principles in Modal Logic \\ and the Axiom of Choice}

\date{}

\author{Rodrigo Nicolau Almeida}
\address{\parbox{\linewidth}{Institute of Logic, Language and Computation\\
University of Amsterdam\\
Amsterdam, 1098XH\\
The Netherlands\\
 }}
\email{r.dacruzsilvapinadealmeida@uva.nl}

\author{Guram Bezhanishvili}
\address{\parbox{\linewidth}{Department of Mathematical Sciences\\
New Mexico State University\\
Las Cruces, NM 88003\\
USA\\
 }
}
\email{guram@nmsu.edu}

\begin{abstract}
We investigate the set-theoretic strength 
of several maximality principles that play an important role in the study of
modal and intuitionistic logics.
We focus on the well-known Fine 
and Esakia maximality principles, present two formulations of each,  
and show that the stronger formulations
are equivalent to the Axiom of Choice (AC), while the weaker ones
to the Boolean Prime Ideal Theorem (BPI). 
\end{abstract}

\subjclass[2020]{03E25; 03B20; 03B44; 03B45; 06D20; 06D22; 06D50; 06E25; 06E15; 18F70}

\keywords{Axiom of Choice; Maximal Ideal Theorem; Prime Ideal Theorem; distributive lattice; Heyting algebra; bi-Heyting algebra; locale; modal algebra; temporal algebra; Priestley duality; Esakia duality; J\'{o}nsson-Tarski duality}

\maketitle

\tableofcontents

\section{Introduction}

In the study of modal and intuitionistic logics --- like in so many fields of mathematics --- it is common to appeal to ``maximality principles," which ensure that certain maximal elements exist. For example, Fine's completeness results in modal logic 
make use of the existence of ``maximal points" in definable subsets of 
canonical models \cite{Fine1974,Fine1985}, while Esakia's analysis of duals of Heyting algebras (now called \textit{Esakia spaces}) proves the existence of maximal points in nonempty closed subsets of such spaces (see, e.g., \cite[Thm.~3.2.1]{Esakiach2019HeyAlg}). The pervasiveness of such principles is often accompanied by the fact that their proofs invoke AC.

It has long been known that some such maximality principles are indeed equivalent to AC, 
while others to its strictly weaker version in the form of BPI
(see, for example, Ern\'{e}'s article \cite{Erné2006} on the subject). In turn, BPI is equivalent to the Prime Ideal Theorem for distributive lattices (see, e.g., \cite[p.~117]{Ern2000} and the references therein): 
\begin{align*}
    \tag{PIT}
    &\text{ For each filter $F$ and ideal $I$ with $F\cap I=\emptyset$, there is}\\
    &\text{a prime ideal $P\supseteq I$ such that $F\cap P=\emptyset$.}
\end{align*}

Following \cite[p.~3]{johnstone1982stone},
all lattices considered in this paper are assumed to be bounded, i.e.~having a top and bottom.
Of particular importance to us is the following 
principle considered by Scott \cite{scotttheorem} (see also Klimovsky \cite{klimovskyexistenceofmaximalfilters} and Mr\'{o}wka \cite{Mrowka1959}), which we will refer to as
\textit{Scott's principle}:
\begin{equation}
        \tag{ScMax} \text{ Each nontrivial distributive lattice has a maximal ideal.}
    \end{equation}


If we work with filters instead of ideals, we obtain a similar principle:
\begin{equation}
        \tag{DLMax} \text{ Each nontrivial distributive lattice has a maximal filter.}
    \end{equation}



Since the order-dual of each distributive lattice ($D,\le)$ is the distributive lattice $(D,\ge)$ whose maximal filters correspond to maximal ideals of $(D,\le)$, we obtain that AC is also equivalent to DLMax. We thus have:

\begin{theorem}\label{Klimovskystheorem}
    The following are equivalent:
    \begin{enumerate}
        \item \textnormal{ScMax};
        \item \textnormal{DLMax};
        \item \textnormal{AC}. 
    \end{enumerate}
\end{theorem}

That the existence of maximal ideals in distributive lattices implies choice is in contrast with what happens in the setting of Heyting algebras. Indeed, the statement that ``each nontrivial Heyting algebra has a maximal ideal" remains equivalent to AC \cite[p.~36]{johnstone1982stone}, but since the order-dual of a Heyting algebra may not be a Heyting algebra, the statement that ``each nontrivial Heyting algebra has a maximal filter" turns out to be equivalent to BPI \cite[p.~79]{johnstone1982stone}.

As was shown by Herrlich \cite{Herrlich2003},
AC remains equivalent to the existence of maximal ideals in the smaller class of those Heyting algebras that are isomorphic to 
the algebras $\mathcal O(X)$ of open sets of topological spaces $X$. The complete Heyting algebras and the ones of the form $\mathcal O(X)$ 
are known as {\em locales} and {\em spatial locales}, and are primary objects of study in pointfree topology \cite{johnstone1982stone,PP2012}. This leads to the following maximality principles:

\begin{enumerate}
    \item[] The \textit{Heyting Maximality Principle}:
    \begin{equation}
        \tag{HMax} \text{ Each nontrivial Heyting algebra has a maximal ideal.}
    \end{equation}
    \item[] The \textit{Locale Maximality Principle}:
    \begin{equation}
        \tag{LMax} \text{ Each nontrivial locale has a maximal ideal.}
    \end{equation}
    \item[] The \textit{Spatial Locale Maximality Principle}:
    \begin{equation}
        \tag{SLMax} \text{ Each nontrivial spatial locale has a maximal ideal.}
    \end{equation}
\end{enumerate}

We again have that each of these is equivlant to AC:

\begin{theorem} 
    The following principles are equivalent:
    \begin{enumerate}
        \item \textnormal{HMax};
        \item \textnormal{LMax};
        \item \textnormal{SLMax};
        \item \textnormal{AC}.
    \end{enumerate}
\end{theorem}

As was pointed out by Banaschewski (see \cite[p. 324]{Herrlich2003}), LMax can be strengthened by reducing the class of locales to compact locales. We can further strengthen this by reducing it to even coherent locales --- an important subclass of compact locales closely related to distributive lattices. Indeed,  
a locale is {\em coherent} provided it is isomorphic to the locale of all ideals of a distributive lattice 
\cite[p.~64]{johnstone1982stone}. We then have the following coherent locale maximality principle:



    \begin{equation}
        \tag{CLMax}
        \text{Each nontrivial coherent locale has a maximal ideal.}
   \end{equation}

\begin{theorem} 
    \textnormal{CLMax} is equivalent to \textnormal{AC}.
\end{theorem}
\begin{proof}
    By Theorem~\ref{Klimovskystheorem}, it is sufficient to show that CLMax is equivalent to ScMax. Clearly ScMax implies CLMax. For the other implication, let $D$ be a nontrivial distributive lattice. We show that it has a maximal ideal. Consider the coherent locale $\mathsf{Id}(D)$ of all ideals of $D$. By CLMax, there is a maximal ideal  $M$ of $\mathsf{Id}(D)$. Let $e:D\to \mathsf{Id}(D)$ be the embedding given by $e(a)={\downarrow}a$ for each $a\in D$.
    We show that 
        $e^{-1}[M]$ is a maximal ideal in $D$.
        Let
    $e^{-1}[M]\subset N$ where $N\in{\sf Id}(D)$. If $a\in N\setminus e^{-1}[M]$, then $e(a)\notin M$, so by maximality of $M$, 
    there is $I\in M$ such that ${\downarrow}a \vee I=D$. 
    This means that $a\vee b=1$ for some $b\in I$ (see, e.g., \cite[p.58]{balbes1974distributive}). Since $b \in I$, we have ${\downarrow}b \subseteq I \in M$, so ${\downarrow}b \in M$, and hence $b \in e^{-1}[M] \subseteq N$. Thus, $1 = a\vee b \in N$, yielding that $N=D$. Consequently, $e^{-1}[M]$ is maximal.
    Thus, ScMax holds.
\end{proof}


The above maximality principles are {\em global} in that they state the existence of maximal ideals (or filters) 
in a given algebra. We will also consider {\em relativized} maximality principles, which state the existence of maximal filters/ideals not only in the algebra, but in a class of algebras related to it, such as homomorphic images of some subreduct. Indeed, it is the relativized maximality principles that are directly related to Esakia's maximality principle. For example, as we will see, the algebraic form of 
this principle states that ``every nontrivial distributive lattice that is a homomorphic image of the lattice reduct of a Heyting algebra has a maximal filter." 
This principle is equivalent to AC, which is in contrast with the corresponding global principle that each nontrivial Heyting algebra has a maximal filter. Indeed, as we already saw, the latter is equivalent to BPI. 


A similar situation can be observed in the setting of 
BAOs (Boolean algebras with operators) of J\'onsson and Tarski \cite{jonnsontarskibaos}. Here we are mainly, but not exclusively, interested in those BAOs that serve as algebraic models of 
the basic modal logic $\mathsf{K}$, as well as their subclasses, with special emphasis on $\mathsf{K4}$-algebras and $\mathsf{S4}$-algebras. The latter were introduced under the name of \textit{closure algebras} by McKinsey and Tarski \cite{mckinseytarskithealgebraoftopology}.
In this setting, the relevant notion of ``maximality" is replaced by ``quasi-maximality" since ultrafilters admit a natural quasi-order. These algebras are closely
related to Heyting algebras, which is at the heart of faithfulness of the G\"{o}del translation of intuitionistic logic into $\mathsf{S4}$ (see, e.g., \cite{rasiowa1968mathematics}). 
For these and related algebras, we will study various versions of 
Fine's maximality principle. 



Our main contribution is to clarify the 
strength of
these naturally occurring principles in modal logic and related areas.
We will show that many of these 
are equivalent to 
AC, 
while others either follow from or are equivalent to BPI.
This can be done by either  
emphasizing the algebraic side or 
by working with the dual spaces of the corresponding algebras, which brings the order-topological side to the fore. 
In the article, we emphasize both perspectives. 

For this purpose, we work with generalizations of Stone duality to distributive lattices, Heyting algebras, and BAOs. 
Such generalizations were initiated by Stone himself, when he generalized his celebrated duality for Boolean algebras to
distributive lattices \cite{Stone1938}. This has resulted in what are now called {\em spectral spaces} \cite{Hochster1969PrimeIS,Dickmann2019}, which 
play a central role 
in algebraic geometry as the spectra of commutative rings \cite{Hartshorne1977}. 
Our focus, however, will be on ordered topological spaces and Priestley duality for distributive lattices \cite{Priestley1970,Priestley1972}
as well as on Stone
spaces equipped with continuous relations and J\'onsson-Tarski 
and Esakia dualities \cite{jonnsontarskibaos,esakiatopologicalkripke}.

We emphasize that Stone duality and its generalizations rely on BPI (see \cite[p.~119]{johnstone1982stone} and Remark \ref{Equivalence of PIT and Priestley duality} below). Thus, in the topological formulations of our results we will always assume BPI, while the algebraic formulations are within $\mathsf{ZF}$. To aid the reader, the statements explicitly mention whether they rely on BPI or are provable within ZF.

\section{Esakia's Maximality Principle}\label{Relativised maximality principles}

\subsection{Distributive lattices and Priestley spaces}


Let $X$ be a topological space. We recall that $U\subseteq X$ is {\em clopen} if $U$ is both closed and open, that $X$ is {\em zero-dimensional} if clopen sets form a basis for the topology, and that $X$ is a {\em Stone space} if $X$ is compact, Hausdorff, and zero-dimensional.

It is a  celebrated result of Stone \cite{Stone1936} that each Boolean algebra is represented as the Boolean algebra $\mathsf{Clop}(X)$ of clopen subsets of a Stone space $X$, which is unique up to homeomorphism. The Stone space $X$ is constructed as the space of ultrafilters of $B$ and the topology on $X$ is generated by the basis $\{ \sigma(a) : a \in B \}$, where $\sigma : B \to \wp(X)$ is the Stone map $\sigma(a) = \{ x \in X : a \in x \}$. The proof that $\sigma$ is an embedding requires BPI, thus Stone's result depends on BPI:

\begin{theorem}[Stone representation] {\em [ZF+BPI]}
    \label{Stone duality}
    Each Boolean algebra is isomorphic to the Boolean algebra of clopen sets of a Stone space.
\end{theorem}

Priestley \cite{Priestley1970} generalized this result by representing distributive lattices as lattices $\mathsf{ClopUp}(X)$ of clopen upsets of the ordered Stone spaces which became known as Priestley spaces.

We recall that an \textit{ordered topological space} is a topological space $X$ equipped with a partial order $\le$. 
For $S\subseteq X$, we let \[
{\uparrow}S=\{x \in X : \exists y\in S \mbox{ with } x\geq y\}\] and
\[{\downarrow}S=\{x \in X : \exists y\in S \mbox{ with } x\leq y.\}\]
Then $S\subseteq X$ is an \textit{upset} if $S={\uparrow}S$, and a \textit{downset} if $S={\downarrow}S$.





\begin{definition}
    Let $X$
    be an ordered topological space. 
    \begin{enumerate}
        \item We say that $X$ satisfies the \textit{Priestley separation axiom} if whenever $x\nleq y$, there is a clopen upset $U$ such that $x\in U$ and $y\notin U$.
        \item We call $X$ a \textit{Priestley space} if $X$ is compact and satisfies the Priestley separation axiom.
    \end{enumerate}
\end{definition}

It is well known (see, e.g., \cite[pp.~257--258]{Davey2002-lr}) that every Priestley space is a Stone space, and that $\leq$ is a closed subset of $X\times X$. In particular, 
${\uparrow}x$ and ${\downarrow}x$ are closed for each $x\in X$.
Priestley spaces play the same role with respect to distributive lattices as Stone spaces with respect to Boolean algebras:

\begin{theorem}[Priestley representation] {\em [ZF+BPI]}
    \label{Priestley duality Theorem}
    Each distributive lattice is isomorphic to the lattice of clopen upsets of a Priestley space.
\end{theorem}

The use of BPI in the Priestley representation theorem is twofold. 
Recall that the Priestley space $X_D$ of a distributive lattice $D$ is constructed as the poset of prime filters of $D$, ordered by inclusion, and the topology has \[\{ \sigma(a) : a \in D \} \cup \{ \sigma(a)^c : a \in D \}\] as a subbasis. 
On the one hand, we require PIT, which is equivalent to BPI (see the introduction), 
to prove that $\sigma$ is an order-embedding. On the other hand, BPI is needed to 
prove that the topology is compact.



\begin{remark}\label{Equivalence of PIT and Priestley duality}\
\begin{enumerate}
    \item As was pointed out in \cite
    {johnstone1982stone}, not only we need BPI to prove Priestley representation,  but the two are equivalent. Indeed, if the Priestley space of each distributive lattice is compact, then so is the Stone space of each Boolean algebra, and the latter implies BPI (see
    \cite[p.~119]{johnstone1982stone}).
    \item The equivalence of the Priestley representation theorem, PIT, and BPI can likewise be seen through the use of \textit{Boolean envelopes} (see \cite[p.~97]{balbes1974distributive}). Indeed, to prove PIT from BPI, it suffices to look at the Boolean envelope $B(D)$ of a distributive lattice $D$; a prime ideal in this Boolean algebra 
    restricts to a prime ideal in $D$.
\end{enumerate}
\end{remark}



For a distributive lattice $D$, there is a well-known dual isomorphism between $X_D$ and the poset $Y_D$ of prime ideals (also ordered by inclusion), obtained by sending each
$x \in X_D$ to 
$D \,\setminus\, x \in Y_D$.
Under this dual isomorphism, maximal ideals of $D$ correspond to minimal prime filters. This allows us to reformulate Theorem \ref{Klimovskystheorem}
in order-topological terms. For a poset $X$, let $\max X$ and $\min X$ be the sets of maximal and minimal points of $X$, respectively. Consider the following maximality and minimality principles for Priestley spaces:


\begin{definition}\ \label{def: PMax and PMin}
\begin{enumerate}
    \item The \textit{Priestley Maximality Principle}:
    \begin{equation*}
    \tag{PMax}
    \text{ If $X$ is a nonempty Priestley space, then $\max X\neq \emptyset$.}
\end{equation*}
    \item The \textit{Priestley Minimality Principle}:
    \begin{equation*}
    \tag{PMin}
    \text{ If $X$ is a nonempty Priestley space, then $\min X\neq \emptyset$.}
\end{equation*}
\end{enumerate}
\end{definition}



By Priestley duality and Theorem \ref{Klimovskystheorem}, PMax and PMin are also equivalent to AC. At the core of this lies the well-known fact that for each distributive lattice ${(D,\leq)}$, its order-dual $(D,\geq)$ is also a distributive lattice, and hence the class $\sf DLat$ of all distributive lattices is closed under taking order-duals.
The maximum of the Priestley space of $(D,\leq)$ is in a one-to-one correspondence with the minimum of the Priestley space of $(D,\geq)$, thus yielding that (PMax) and (PMin) are equivalent.

\begin{remark}
    In light of the above discussion, each of the maximality principles is equivalent to the corresponding minimality principle: instead of considering the existence of maximal ideals, one considers the existence of minimal prime filters.
\end{remark}


\subsection{Heyting algebras and Esakia spaces}

\begin{definition}\
A \textit{Heyting algebra} is a 
   lattice $H$ with an additional binary operation $\to:H^2\to H$ such that, 
   for all $a,b,c\in H$,
   \begin{equation*}
       a\wedge c\leq b\iff c\leq a\rightarrow b.
   \end{equation*}
\end{definition}



It is well known that each Heyting algebra is a distributive lattice, and that 
Priestley representation 
specializes to Esakia representation of Heyting algebras \cite{esakiatopologicalkripke}:

\begin{definition}
    We say that a Priestley space $X$ is an \textit{Esakia space} if whenever $U$ is a clopen subset, then ${\downarrow}U$ is clopen.
\end{definition}





%


\begin{theorem} [ZF+BPI]
$ $
{\em (Esakia representation)} 
    Each Heyting algebra is isomorphic to the Heyting algebra of clopen upsets of an Esakia space.
\end{theorem}

As noted in the introduction, unlike distributive lattices, the existence of maximal filters in Heyting algebras follows from (and is equivalent to) BPI. Indeed, 
for each Heyting algebra $H$, there is a bijection between the maximal filters of $H$ and the maximal filters of its {\em Booleanization} 
\[
\mathfrak B(H):=\{a\in H : a=\neg\neg a\},
\]
which is a Boolean algebra (see, e.g., \cite[p.~75]{johnstone1982stone}). 
In light of Esakia duality, this can be stated as follows:

\begin{proposition}[ZF+BPI]\label{Dual version of Johnstones theorem}
    If $X$ is an Esakia space, then $\max X\neq \emptyset$.
\end{proposition}

The reason for this is that, unlike distributive lattices, the order-dual of a Heyting algebra may not be a Heyting algebra. 
Thus, the existence of maximal ideals in Heyting algebras does not follow from the above. We will now see that the picture changes when one focuses instead on \textit{relativized} principles, where not only the algebra, but a whole class of its homomorphic images, are asked to contain maximal ideals.

\subsection{Relativized maximality principles for Esakia spaces} 

The first relativized maximality principle we will consider was introduced by Esakia 
in his studies of 
Heyting algebras and their dual spaces. Instead of requiring that $\max X$ is nonempty for each nonempty Esakia space $X$,
we \textit{relativize} this requirement to all nonempty closed subsets of $X$:

\begin{definition} \cite{Esakiach2019HeyAlg}
    {\em Esakia's Maximality Principle}:
\begin{equation}
    \tag{REMax}
    \text{ For each Esakia space $X$ and a nonempty closed subset $C\subseteq X$, $\max C\neq \emptyset$.}
\end{equation}    
\end{definition}

Esakia proved that AC implies REMax (see \cite[Thm.~3.2.1]{Esakiach2019HeyAlg}). To prove the converse, 
we make use of some facts about \textit{free distributive lattices} (see, e.g., \cite[Sec.~V.3]{balbes1974distributive}). For a set $X$, let $F_{\sf DLat}(X)$ be the free distributve lattice generated by $X$. We denote by $\mathbf{2}^{X}$ the ordered topological space which is the product of $X$ many copies of the two element chain $\mathbf{2}=\{0,1\}$, with the product topology and order.

\begin{proposition}(ZF+BPI)\label{Product of Esakia space}
    The Priestley dual of $F_{\sf DLat}(X)$ is isomorphic to $\mathbf{2}^{X}$, which is an Esakia space.
\end{proposition}
\begin{proof}
    The first statement can be found in \cite[Thm.~202]{Gratzer2011-qy}, and the second follows from 
    the fact that a product of Esakia spaces is an Esakia space \cite[Prop.A.8.1]{Esakiach2019HeyAlg}.
\end{proof}

We say that two Priestley spaces $X$ and $Y$ are {\em isomorphic} if there is a homeomorphism $f:X\to Y$ which in addition is an order-isomorphism. The following lemma is well known  
(see, e.g., \cite[Thm.~11.31(ii)]{Davey2002-lr}). 

\begin{lemma}\label{Duality for surjective homomorphisms} (ZF+BPI)
    Let $D,D'$ be distributive lattices and $f:D\to D'$ an onto homomorphism. Then $X_{D'}$ is isomorphic to a closed subspace of $X_{D}$.
\end{lemma}

\begin{theorem}[ZF+BPI]\label{Equivalence of Esakias Maximality with the Axiom of Choice}
    \textnormal{REMax} is equivalent to \textnormal{AC}.
\end{theorem}
\begin{proof}
As we already pointed out, the implication AC $\Rightarrow$ REMax is proved in \cite[Thm.~3.2.1]{Esakiach2019HeyAlg}.
We show the implication REMax $\Rightarrow$ AC. It suffices to prove PMax.
Let $X$ be a nonempty 
Priestley space, and let $D$ be the distributive lattice of clopen upsets of $X$.
Then $D$ is a homomorphic image of $F_{\sf Dist}(Z)$ for some set $Z$. By Proposition \ref{Product of Esakia space} and Lemma \ref{Duality for surjective homomorphisms}, $X$ is isomorphic to a nonempty closed subspace of $\mathbf{2}^{Z}$, which 
is an Esakia space. Thus, by REMax,
$\max X\ne\varnothing$. 
\end{proof}


The above equivalence relies on BPI. We next provide a more algebraic formulation of Esakia's Maximality Principle and prove its equivalence with AC without invoking BPI.

\begin{definition}
    Let $H$ be a Heyting algebra. We say that a  distributive lattice $D$ is a \textit{lattice homomorphic image} of $H$ if there is an onto lattice homomorphism $f:H\to D$. We write $\mathbb{H}_{\mathsf{Dist}}(H)$ for the class of all lattice homomorphic images of~$H$.
\end{definition}

Note that 
even if $D\in \mathbb{H}_{\mathsf{Dist}}(H)$ happens to be a Heyting algebra, the lattice homomorphism $h : H \to D$ 
need not be a Heyting homomorphism. 

\begin{definition}
    \textit{Relativized Heyting Maximality Principle}:
    \begin{align*}
        \tag{RHMax}
        &\text{For each Heyting algebra $H$ and a lattice homomorphic}\\ &\text{image $D$ of $H$, $D$ has a maximal filter.}
    \end{align*}
\end{definition}

Let $H$ be a Heyting algebra and $X_{H}$ its dual space. By Priestley duality, lattice homomorphic images of $H$ correspond to closed subsets of $X_H$ (see Lemma \ref{Duality for surjective homomorphisms}). Thus, in ZF+BPI, RHMax is equivalent to REMax, and hence to AC by Theorem \ref{Equivalence of Esakias Maximality with the Axiom of Choice}. We next show that the equivalence of RHMax and AC can be proved without invoking BPI. For this we prove that $F_{\sf Dist}(X)$ is a Heyting algebra. Of course, this follows from Proposition \ref{Product of Esakia space}, which uses BPI. We prove this constructively.

\begin{lemma}\label{Dual is a Heyting algebra for each X algebraically}
$F_{\sf Dist}(X)$ is a Heyting algebra for each set $X$.
\end{lemma}
\begin{proof}
Let $P,Q \in F_{\mathsf{Dist}}(X)$. We show that the implication $P\to Q$ exists in $F_{\mathsf{Dist}}(X)$.
By \cite[p.~89, Thm.~7]{balbes1974distributive}, 
    $P=\bigvee_{i=1}^{n}S_{i}$, where each $S_{i}$ is join-irreducible, and $Q=\bigwedge_{j=1}^{m}T_{i}$, where each $T_{j}$ is meet-irreducible. 
    It is sufficient 
    to show that each $S_i \to T_{j}$ exists in $F_{\mathsf{Dist}}(X)$ since then 
    \begin{eqnarray*}
        P\rightarrow Q &= \left(\bigvee_{i=1}^{n}S_{i}\right) \rightarrow \left(\bigwedge_{j=1}^{m}T_{i}\right) \\
        &= \bigwedge_{i=1}^{n} \left(S_{i} \rightarrow \bigwedge_{j=1}^{m}T_{i}\right) \\
        &= \bigwedge_{i=1}^{n}\bigwedge_{j=1}^{m}\left(S_{i}\rightarrow T_{j}\right).
    \end{eqnarray*}
     
We show that
\[
S_{i}\rightarrow T_{j}=
\begin{cases}
    1 \text{ if } S_{i}\leq T_{j};\\
    T_j \text{ otherwise. }
\end{cases}
\]
To see this, assume that $S_{i}\nleq T_{j}$. Let $R \in F_{\mathsf{Dist}}(X)$ with $S_{i}\wedge R\leq T_{j}$. Since $T_{j}$ is meet-irreducible, our assumption yields that 
$R\leq T_{j}$. 
Thus, $S_i \to T_{j} = T_j$, concluding the proof. 
\end{proof}

    \begin{theorem} [ZF]
        \textnormal{RHMax} is equivalent to \textnormal{AC}.
    \end{theorem}
    \begin{proof}
        It is clear that AC implies RHMax since it implies the distributive lattice maximality principle.  For the other implication, let $D$ be 
        a nontrivial distributive lattice. Then 
        $D$ is a homomorphic image of $F_{\sf Dist}(Z)$ for some set $Z$. By Lemma~ \ref{Dual is a Heyting algebra for each X algebraically}, $F_{\sf Dist}(Z)$ is a Heyting algebra, which is non-trivial since $D$ is a non-trivial distributive lattice. Thus, by RHMax,
        $D$ has a maximal filter. This shows DLMax, which entails AC by 
        Theorem \ref{Klimovskystheorem}. 
\end{proof}

\section{Fine's Maximality Principle}\label{Relativised maximality in modal logic}

\subsection{Maximality principles in $\mathsf{K4}$ and $\mathsf{S4}$}\label{Fines Maximality Principle}

To formulate Fine's Maximality Principle, we recall some basic facts from modal logic \cite{Chagrov1997-cr,Blackburn2002-fd}.

\begin{definition}
    A \textit{modal algebra} is a Boolean algebra $B$
    equipped with a unary function $\Box:B\to B$ preserving finite meets.
\end{definition}

Modal algebras are a special case of Boolean algebras with operators (BAOs) of J\'{o}nnson and Tarski \cite{jonnsontarskibaos}.
In this section we will be focusing on  
modal algebras satisfying additional conditions (see, e.g., \cite[p.~214]{Chagrov1997-cr}):

\begin{definition}
    A modal algebra $(B,\Box)$ is \begin{enumerate}
        \item a {\em $\mathsf{K4}$-algebra} if $\Box a\leq \Box\Box a$ for each $a\in B$;
        \item an {\em $\mathsf{S4}$-algebra} if additionally it satisfies $\Box a\leq a$ for each $a\in B$.
    \end{enumerate}
\end{definition}



%




J\'onnson and Tarski \cite{jonnsontarskibaos} extended Stone representation of Boolean algebras to BAOs. As a special case, this yields the following representation of modal algebras. We recall that a binary relation $R$ on a Stone space $X$ is {\em continuous} if 
\begin{samepage}
\begin{itemize}
    \item $R[x] := \{y\in X : xRy\}$ is closed for each $x\in X$, and
    \item $R^{-1}[U] := \{x\in X : \exists y\in U$ with $xRy\}$ is clopen for each clopen $U \subseteq X$. 
\end{itemize}
\end{samepage}
 A {\em modal space} is a Stone space $X$ equipped with a continuous relation $R$. If $(X,R)$ is a modal space, then $({\sf Clop}(X),\Box_R)$ is a modal algebra, where
\begin{equation*}
    \Box_{R}U=\{x\in X : R[x] \subseteq U\}
\end{equation*}
for each $U\in{\sf Clop}(X)$. Moreover, we have: 

\begin{theorem}[J\'onsson-Tarski representation] {\em [ZF+BPI]} \label{Representation of Modal algebras}
    Each modal algebra $(B,\Box)$ is represented as the modal algebra $({\sf Clop}(X),\Box_R)$ for a modal space $(X,R)$.
\end{theorem}

The representation space of $(B,\Box)$ is the pair $(X_B,R_\Box)$, where $X_B$ is the Stone space of $B$ and 
\begin{equation*}
    xR_{\Box}y \iff (\forall a\in B) \, (\Box a\in x \rightarrow a\in y).
\end{equation*}





A modal space $(X,R)$ is 
\begin{itemize}
    \item a \textit{$\sf K4$-space} if $R$ is a transitive relation, and 
    \item an \textit{$\sf S4$-space} if $R$ is in addition a reflexive relation, so $R$ is a quasi-order. 
\end{itemize}
Observe that Esakia spaces are exactly the partially ordered $\sf S4$-spaces. 
The following 
goes back to J\'{o}nnson and Tarski \cite{jonnsontarskibaos}:

\begin{theorem} [ZF+BPI]
    A modal algebra $(B,\Box)$ is
\begin{enumerate}
        \item a $\mathsf{K4}$-algebra iff $(X_{B},R_\Box)$ is a $\mathsf{K4}$-space;
        \item an $\mathsf{S4}$-algebra iff $(X_{B},R_\Box)$ is an $\mathsf{S4}$-space.
    \end{enumerate}
\end{theorem}

For a modal space $(X,R)$ and $C\subseteq X$, let $R_{C}$ be the restriction of $R$ to $C$; that is, $R_{C}=R\cap (C\times C)$. The following is well known (see, e.g.,
\cite[Thm.~3.2.6]{Esakiach2019HeyAlg}):

\begin{theorem} \label{Inheritancetheoremfork4algebras}
    Let $(X,R)$ be a modal space. If $C$ is a clopen subset of $X$, then $(C,R_C)$ is a modal space. Moreover, if $(X,R)$ is a $\mathsf{K4}$-space, then so is $(C,R_C)$; and if $(X,R)$ is an $\mathsf{S4}$-space, then so is $(C,R_C)$.
\end{theorem}

There is a close connection between 
$\mathsf{K4}$-algebras and $\mathsf{S4}$-algebras. For a $\mathsf{K4}$-algebra $(B,\Box)$, define $\Box^{+}:B\to B$ by  $\Box^{+}a=a\wedge\Box a$ for each $a\in B$. Then $(B,\Box^{+})$ is an $\mathsf{S4}$-algebra \cite[Appendix I]{mckinseytarskithealgebraoftopology}.
Dually,
if $(X,R)$ is a $\mathsf{K4}$-space, then $(X,R^{+})$ is an $\mathsf{S4}$-space, where $R^{+}=R\cup \{(x,x) : x\in X\}$ is the {\em reflexivization} of $R$ (see, e.g., \cite[Lem.~3.6]{bezhanishvilighilardijibladze}).



Let $(X,R)$ be a quasi-ordered set. We recall (see, e.g., \cite[Def.~1.4.9]{Esakiach2019HeyAlg}) that a point $x\in X$ is \textit{quasi-maximal} if $xRy$ implies $yRx$. Let $\mathsf{qmax} X$ be the set of quasi-maximal points of $X$. 
We will mainly be interested in quasi-maximal points of $\mathsf{S4}$-spaces and $\mathsf{K4}$-spaces. For the latter, the definition is the same. Observe that 
$\mathsf{qmax}(X,R)=\mathsf{qmax}(X,R^{+})$.

We are ready to formulate Fine's Maximality Principle. We recall that to prove his completeness results, Fine \cite{Fine1974,Fine1985} worked with maximal points in definable subsets of canonical models of the modal logic $\sf K4$ and its extensions. In the language of $\sf K4$-spaces, definable subsets correspond to clopens. Thus, Fine's Maximality Principle can be formulated as follows.

\begin{definition}
The \textit{Fine Maximality Principle}: 
    \begin{align*}
    \tag{FMax}
        &\text{ If $(X,R)$ is a $\sf K4$-space and $C\subseteq X$ is a}\\
        &\text{ nonempty clopen, then $\mathsf{qmax} C\neq \emptyset$.}
    \end{align*}
\end{definition}


We prove that FMAx is a consequence of BPI. We then provide a stronger formulation of Fine's Maximality Principle for nonempty closed subsets of $\sf K4$-spaces, and show that it is equivalent to AC. For this, 
we exploit the intimate connection between $\mathsf{S4}$-algebras and Heyting algebras observed by 
McKinsey and Tarski \cite{mckinseytarskiclosedelements} (see also \cite{rasiowa1968mathematics,Esakiach2019HeyAlg,Chagrov1997-cr}). If $(B,\Box)$ is an $\mathsf{S4}$-algebra, then the $\Box$-fixpoints 
\begin{equation*}
    B_{\Box}=\{a\in B : a=\Box a\}
\end{equation*}
form a Heyting algebra (see, e.g., \cite[Prop.~2.2.4]{Esakiach2019HeyAlg}).

The Esakia space of $B_\Box$ has the following convenient description. 
Let $(X,R)$ be an $\sf S4$-space. Define the equivalence relation of being in the same cluster:
\begin{equation*}
    x\sim y \iff xRy \text{ and } yRx.
\end{equation*}
Let $\rho(X)\coloneqq X/{\sim}$ be the 
quotient space and $\rho:X\to X/\sim$ the quotient map. Define 
\[
[x]\leq [y] \iff xRy.
\]
Then $\leq$ is well defined, and $(\rho(X),\leq)$ is an Esakia space. In fact, if $(X,R)$ is the dual space of an $\sf S4$-algebra $(B,\Box)$, then $(\rho(X),\leq)$ is (isomorphic to) the Esakia space of the Heyting algebra $B_\Box$ (see \cite[pp.~65--66]{Esakiach2019HeyAlg}).



\begin{theorem}[ZF+BPI]\ \label{FMax and existence of quasimaximal points}
    \begin{enumerate}
        \item Every nonempty $\mathsf{S4}$-space
        has a quasi-maximal point.
        \item Every nonempty $\mathsf{K4}$-space
        has a quasi-maximal point.
        \item \textnormal{FMax} holds.
    \end{enumerate}
\end{theorem}
\begin{proof}
    (1) 
    Let $(X,R)$ be a nonempty $\sf S4$-space. Then $(\rho(X),\leq)$ is a nonempty Esakia space. 
    Since we assume BPI, Proposition \ref{Dual version of Johnstones theorem} applies, by which
    $\rho(X)$ has a maximal point, say $[x]$.
    Then $x$ must be a quasi-maximal point of $X$ since if $xRy$, then $[x]\leq [y]$, so $[x]=[y]$, and hence $yRx$.

    (2) Let $(X,R)$ be a nonempty $\sf K4$-space. 
    Then $(X,R^{+})$ is a nonempty $\sf S4$-space. By (1), there is a quasi-maximal point $x$ in $(X,R^+)$. 
    But then $x$ is quasi-maximal in $(X,R)$. 

    (3) Let $(X,R)$ be a $\sf K4$-space and $C\subseteq X$ a nonempty clopen. By Theorem~\ref{Inheritancetheoremfork4algebras}, $(C,R_{C})$ is a $\sf K4$-space. By (2) we then have that $(C,R_{C})$ has a quasi-maximal point, proving FMax.
\end{proof}

Next we can strengthen Fine's Maximality Principle:

\begin{definition}\
\begin{enumerate}
    \item \textit{Fine's Strong Maximality Principle}:
    \begin{align*}
    \tag{FSMax}
        &\text{ If $(X,R)$ is a $\sf K4$-space and $C\subseteq X$ is a nonempty }\\
        &\text{closed subset, then $\mathsf{qmax} C\neq \emptyset$.}
    \end{align*}
    \item \textit{Fine's Strong Maximality Principle for $\mathsf{S4}$}:
    \begin{align*}
    \tag{$\mathrm{FSMax^{S4}}$}
        &\text{ If $(X,R)$ is an $\sf S4$-space and $C\subseteq X$ is a nonempty }\\
        &\text{closed subset, then $\mathsf{qmax} C\neq \emptyset$.}
    \end{align*}\end{enumerate}
\end{definition}

\begin{theorem}[ZF+BPI]\label{Equivalence of Strong Maximality and Axiom of Choice}
The following principles are equivalent:
\begin{enumerate}
    \item \textnormal{FSMax};
    \item $\textnormal{FSMax}^{S4}$;
    \item \textnormal{REMax};
    \item \textnormal{AC}. 
\end{enumerate}
\end{theorem}
\begin{proof}
(1)$\Rightarrow$(2) and (2)$\Rightarrow$(3) are immediate since every $\mathsf{S4}$-space is a $\mathsf{K4}$-space, and every Esakia space is an $\mathsf{S4}$-space. By Theorem \ref{Equivalence of Esakias Maximality with the Axiom of Choice}, (3)$\Rightarrow$(4). By \cite[p.~47]{Esakiach2019HeyAlg}, (4)$\Rightarrow$(2). 
We prove (2)$\Rightarrow$(1). Let $(X,R)$ be a $\mathsf{K4}$-space and $C \subseteq X$ a nonempty closed subset. Then $(X,R^{+})$ is an $\mathsf{S4}$-space, and thus $\mathsf{qmax}_{R^{+}}C\neq \emptyset$ by $\mathrm{FSMax^{S4}}$. But then $\mathsf{qmax}_{R}C\neq \emptyset$ since $\mathsf{qmax}_{R}C=\mathsf{qmax}_{R^{+}}C$. 
\end{proof}

Fine's Maximality Principle admits a natural algebraic formulation using \textit{relativizations}. We recall that if $a\in B$, then 
the set
\begin{equation*}
    B\restriction a=\{c : c\leq a\}
\end{equation*}
with the induced order is 
a Boolean algebra, 
and there is an onto Boolean homomorphism $f:B\to B\restriction a$ given by $f(c) = a \wedge c$ \cite[Sec.~II.6]{rasiowa1968mathematics}. We can equip the latter set with a modal operator $\Box_{a}$ by setting $\Box_{a}c=a\wedge \Box (a\to c)$ for each $c\in B\restriction a$ \cite[Sec.~III.2]{rasiowa1968mathematics}.

Note that for a modal algebra $(B,\Box)$, the definition of $R_\Box$ on $X_B$ 
makes sense without BPI (even if $X_B =\varnothing$).
We write $\mathsf{qmax} \, B$ for the set of $R_{\Box}$-quasi-maximal ultrafilters of $B$.



\begin{definition}\
  \begin{itemize}
      \item \textit{Fine's Algebraic Maximality Principle}:
    \begin{align*}
    \tag{FAMax}
        \text{ If } (B,\Box) \text{ is a $\mathsf{K4}$-algebra and $a\in B \setminus \{0\}$, then } \mathsf{qmax} (B\restriction a)\neq \emptyset.
     \end{align*}
     \item \textit{Fine's Algebraic Maximality Principle for $\mathsf{S4}$}:
    \begin{align*}
    \tag{FAMax$^{\mathsf{S4}}$}
        \text{ If } (B,\Box) \text{ is an $\mathsf{S4}$-algebra and $a\in B \setminus \{0\}$, then } \mathsf{qmax} (B\restriction a)\neq \emptyset.
     \end{align*}
  \end{itemize}  
\end{definition}

Let $(B,\Box)$ be a $\mathsf{K4}$-algebra and $(X_{B},R_\Box)$ its dual space. By J\'onsson-Tarski duality, relativizations of $(B,\Box)$ correspond to $(C,R_{C})$, where $C \subseteq X_B$ is clopen. Thus, in ZF+BPI, FAMax is equivalent to FMax. We next show that in ZF, FAMax is equivalent to BPI.


\begin{theorem} [ZF] \label{thm: FAMax} The following are equivalent:
\begin{enumerate}
    \item \textnormal{BPI};
    \item \textnormal{FAMax};
    \item \textnormal{FAMax}$^{\mathsf{S4}}$;
\end{enumerate}
\end{theorem}
\begin{proof}
If we assume (1) then FMax holds by Theorem \ref{FMax and existence of quasimaximal points}(3). Therefore, (2) holds in light of the above paragraph. Certainly (2)$\Rightarrow$(3).

So suppose FAMax$^{\mathsf{S4}}$ holds. It is sufficient to show that every  proper filter in a Boolean algebra can be extended to an ultrafilter.
Let $B$ be a Boolean algebra and $F \subseteq B$ a proper filter. 
Consider the quotient algebra $B/F$,  
and let $p:B\to B/F$ be the quotient map. We can view $B/F$ as a $\mathsf{K4}$-algebra by setting $\Box a=a$ for each $a\in B/F$. Since $B/F= (B/F)\restriction 1$, applying FAMax yields a quasi-maximal ultrafilter $F'\subseteq B/F$. It is easy to see that $p^{-1}[F']$ is then an ultrafilter
in $B$ (see, e.g., \cite[Lem.~IV.3.16] 
{BurrisSankappanavar}), which extends $F$. Thus, BPI holds. 
\end{proof}

It is more difficult to find an algebraic formulation of Fine's Strong Maximality principle.
The key issue is that Theorem \ref{Inheritancetheoremfork4algebras} is not necessarily true if one replaces clopens by closed sets. Namely, the restriction of a continuous relation on a Stone space to a closed subspace may no longer be continuous. Therefore, if $(X,R)$ is a modal space and $C$ is a closed subset of $X$, then $(C,R_C)$ may no longer be a modal space. We give one such example below.

\begin{example}
Consider the space $X$ depicted in Figure \ref{fig:counterexampletoclosedsets}, with the partial order as shown. 
The topology is given by letting each $a_n$ and $b_n$ be isolated, $(a_n)\to a_\omega$, and $(b_n)\to b_\omega$. Thus, $X$ is the two-compactification of a discrete space. It is easy to see that $X$ is a Priestley space in which the downset of each clopen is clopen. Thus, $X$ is an Esakia space, and hence a $\mathsf{K4}$-space. 
Let $C={\downarrow}b_{0}\cup \{a_{\omega}\}$. Then $C$ is a closed subset of $X$, 
but $(C,\leq_C)$ is not a $\mathsf{K4}$-space
because $\{a_{\omega}\}$ is clopen in $C$, but ${\downarrow}\{a_{\omega}\} = \{ a_\omega,b_\omega \}$ is not. 

\begin{figure}[h]
    \centering

    \begin{tikzpicture}
        \node at (0,0) {$\bullet$};
        \node at (-0.3,0) {$a_{0}$};
        \node at (0,-1) {$\bullet$};
        \node at (-0.3,-1) {$a_{1}$};
        \node at (0,-2) {$\bullet$};
        \node at (-0.3,-2) {$a_{2}$};
        \node at (0,-3) {$\vdots$};
        \node at (2,-1) {$\bullet$};
        \node at (2.3,-1) {$b_{0}$};
        \node at (2,-2) {$\bullet$};
        \node at (2.3,-2) {$b_{1}$};
        \node at (2,-3) {$\bullet$};
        \node at (2.3,-3) {$b_{2}$};
        \node at (2,-4) {$\vdots$};
        \node at (0,-4) {$\bullet$};
        \node at (-0.3,-4) {$a_{\omega}$};
        \node at (2,-5) {$\bullet$};
        \node at (2.3,-5) {$b_{\omega}$};

        \draw (0,0) -- (0,-1) -- (0,-2) -- (0,-2.5) -- (0,-2) -- (2,-3) -- (2,-3.5) -- (2,-3) -- (2,-2) -- (0,-1) -- (0,0) -- (2,-1) -- (2,-2);

        \draw (0,-4) -- (2,-5);

    \end{tikzpicture}    \caption{The modal space $(X,R)$ and its closed subspace $C$}
    \label{fig:counterexampletoclosedsets}
\end{figure}
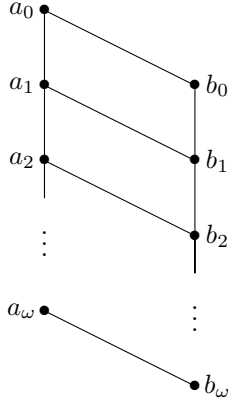
\end{example}

In light of such examples,
a purely algebraic formulation of 
FSMax would be
rather unnatural. 
It is for this reason that in this subsection we have focused 
on the dual 
formulations.




The Fine and Esakia Maximality Principles can further be generalized to weaker modal systems. They can also be studied for temporal modal systems and the closely related bi-Heyting algebras. 
In the rest of the paper we explore both possibilities. 

\subsection{Maximality principles in weaker modal systems}\label{Maximality in weaker modal systems}

We recall (see, e.g., \cite[p.~74]{krachtbook}) that a normal modal logic $\mathsf{L}$ is \textit{n-transitive} provided it proves 
the following principle:
\begin{equation}
\tag{n-tra}
    \Box^{\leq n}p \rightarrow \Box^{\leq n+1}p,
\end{equation}
where $\Box^{\leq n}p$ is the shorthand for $p\wedge \Box p \wedge...\wedge \Box^{n}p$. We say that a normal modal logic $\mathsf{L}$ is \textit{weakly transitive} if it is $n$-transitive for some $n\in \omega$.

These logics were introduced by Rautenberg  \cite{Rautenberg1980}, and have since been studied extensively. 
They are some of the better behaved modal logics since they admit a master modality \cite[p.~74]{krachtbook}, 
satisfy an appropriate form of the deduction theorem, and  
the associated varieties have Equationally Definable Principal Congruences \cite{Blok1982}.

    Let $n\in\omega$. For a modal algebra $(B,\Box)$, define $\Box^n,\Box^{\le n} : B \to B$ by setting $\Box^0 a = \Box a$, $\Box^{n+1}a=\Box\Box^n a$, and $\Box^{\leq n}a=a\wedge\Box a\wedge\cdots\wedge \Box^{n}a$ for each $a\in B$.
    Also, 
    for a binary relation $R$ on a set $X$, set $R^{0}=\{(x,x) : x\in X\}$, 
    $R^{n+1}=R \circ R^n$, and $R^{*}=\bigcup_{n\in \omega}R^{n}$; that is, $R^*$ is the reflexive and transitive closure of $R$. 

\begin{samepage}
\begin{definition}\
    \begin{enumerate}
        \item A modal algebra $(B,\Box)$ is {\em $n$-transitive} if it satisfies $\Box^{\leq n}a\leq \Box^{\leq n+1}a$ for each $a\in B$.
        \item A modal space $(X,R)$ is {\em $n$-transitive} if $R^{n+1}\subseteq R^{n}$.
    \end{enumerate}
\end{definition}
\end{samepage}


The notions of $n$-transitive modal algebras and 
spaces were introduced 
in \cite{Rautenberg1980}, where it was shown that a modal algebra is $n$-transitive iff its dual modal space is $n$-transitive. 
For the next definition, we refer to
\cite{krachtbook}.


\begin{definition}
A modal algebra $(B,\Box)$ is \textit{weakly transitive} if it is $n$-transitive for some $n\in\omega$,
    and a modal space $(X,R)$ is \textit{weakly transitive} if it is $n$-transitive for some $n\in\omega$ (equivalently, 
    $R^{n}=R^{*}$ for some $n\in\omega$).
\end{definition}

Since continuous relations on Stone spaces are closed under composition \cite{Halmos19541956}, we have:

\begin{proposition}\label{Transitive spaces from modal spaces}
    If $(X,R)$ is a modal space, then so is $(X,R^{n})$ for each $n\in\omega$. 
    Moreover, if $R^{*} = R^{n}$ for some $n\in\omega$, then $(X,R^{n})$ is an $\mathsf{S4}$-space. 
\end{proposition}

In general however, if $(X,R)$ is a modal space, $(X,R^{*})$ need not be a modal space since the relation $R^{*}$ need not be continuous (see, e.g., \cite{Venema2004}). 

\begin{proposition}[ZF+BPI] \label{Definition of weakly transitive modal space}
The following are equivalent:
\begin{enumerate}
    \item $(B,\Box)$ is weakly transitive;
    \item $(X_{B},R_{\Box})$ is weakly transitive;
    \item $(X_{B},R^{n}_{\Box})$ is an $\mathsf{S4}$-space for some $n\in\omega$.
\end{enumerate}
\end{proposition}
\begin{proof}
For (1)$\Leftrightarrow$(2) see
\cite{Rautenberg1980}. That (2)$\Rightarrow$(3) follows from Proposition~\ref{Transitive spaces from modal spaces}. 
For (3)$\Rightarrow$(2), if $R_{\Box}^{n}$ is reflexive and transitive, then $R_{\Box}^{n+1}\subseteq R_{\Box}^{n}$. 
\end{proof}



\begin{definition}
    Let $(X,R)$ be a weakly transitive modal space. We call $x\in X$ \textit{eventually quasi-maximal} provided $xR^{*}y$ implies $yR^{*}x$. For a subset $C\subseteq X$, we denote by $\mathsf{eqmax} C$ the set of eventually quasi-maximal points of $C$. 
\end{definition}

We are ready to define the maximality principle, which generalizes Fine's Strong Maximality Principle to weakly transitive modal spaces.

\begin{definition}
    \textit{Weakly Transitive Maximality Principle}:
    \begin{align*}
    \tag{WTMax} &\text{ If $(X,R)$ is a weakly transitive modal space and $C \subseteq X$ }\\
    & \text{ is a nonempty closed subset, then $\mathsf{eqmax}C\neq \emptyset$.}
\end{align*}
\end{definition}


\begin{theorem}[ZF+BPI] \label{Modal Principles equivalent to Axiom of Choice}
    The following are equivalent:
    \begin{enumerate}
    \item \textnormal{WTMax}; 
    \item \textnormal{FSMax}; 
    \item \textnormal{AC}. 
    \end{enumerate}
\end{theorem}
\begin{proof}
    By Theorem \ref{Equivalence of Strong Maximality and Axiom of Choice}, (2) and (3) are equivalent. Note that (1) implies (2) since every $\mathsf{S4}$-space is weakly transitive (indeed, in this case, $R^{*}=R$).  
    We show that (2) implies (1).
    Let $(X,R)$ be a weakly transitive modal space, and let $C\subseteq X$ be a nonempty closed subset. By Proposition \ref{Transitive spaces from modal spaces}, 
    $(X,R^{n})$ is an $\mathsf{S4}$-space for some $n\in\omega$. By FMax, $\mathsf{qmax}_{R^n}C\neq \emptyset$.
    But a point $x$ is quasi-maximal with respect to $(X,R^{n})$ iff it is eventually quasi-maximal with respect to $(X,R)$. Thus, $\mathsf{eqmax}C\neq \emptyset$, concluding the proof.
\end{proof}

\begin{remark}\label{rem: WTMax algebraically}
    As in Section \ref{Fines Maximality Principle}, we could consider the weaker version of WTMax for nonempty clopen subsets. But 
    as in Theorem \ref{FMax and existence of quasimaximal points}, such a principle could be derived from BPI. Moreover, it could be reformulated algebraically to obtain an algebraic maximality principle analogous to FAMax, which is also equivalent to BPI (see Theorem~\ref{thm: FAMax}).
\end{remark}

\subsection{Maximality principles in temporal logic}

The formal study of temporal logic has 
its origins in the work of Arthur Prior \cite{Pri57} and has been studied extensively (see, e.g., \cite{Venema2001-VENTLB,Gabbay1994-ar} for a detailed treatment).

\begin{definition}
    A \textit{tense algebra} is a Boolean algebra $B$ equipped with two unary functions $\Box_{F},\Box_P:B\to B$ 
    such that $(B,\Box_{P})$ and $(B,\Box_{F})$ are modal algebras and the connecting axioms are satisfied, 
    for all $a\in B$:
\begin{equation*}
    a\leq \Box_{F}\Diamond_{P} a \text{ and } a\leq \Box_{P}\Diamond_{F} a,
\end{equation*}
where $\Diamond_{P}a=\neg\Box_{a}\neg a$ and $\Diamond_{F}a=\neg\Box_{a}\neg a$. If additionally $(B,\Box_{F})$ is $\mathsf{S4}$, we say that $(B,\Box_{F},\Box_{P})$ is an {\em $\mathsf{S4}_{t}$-algebra}.
\end{definition}

Intuitively, one reads $\Box_{P}$ as ``always in the past" and $\Box_{F}$ as ``always in the future," and similarly for $\Diamond_{P}$ and $\Diamond_{F}$.

The J\'{o}nnson-Tarski representation of modal algebras extends to tense algebras. For a binary relation $R$, let $R^{\dagger}=\{(y,x) : (x,y)\in R\}$. 
We say that a modal space $(X,R)$ is a \textit{tense space} if the relation $R^{\dagger}$ is also continuous, and we call a tense space an {\em $\mathsf{S4}_{t}$-space} if $R$ is a quasi-order \cite{esakia1975problem}.

If $(X,R)$ is a tense space, then $(\mathsf{Clop}(X),\Box_{R},\Box_{R^{\dagger}})$ is a tense algebra, and we have the following analogue of Theorem \ref{Representation of Modal algebras}:

\begin{theorem}[J\'onsson-Tarski representation] {\em [ZF+BPI]} \label{Representation of tense algebras}
    Each tense algebra $(B,\Box_{F},\Box_{P})$ is represented as the tense algebra $({\sf Clop}(X),\Box_R,\Box_{R^{\dagger}})$ of a tense space $(X,R)$. Moreover, $(B,\Box_{F},\Box_{P})$ is an $\mathsf{S4}_{t}$-algebra iff $R$ is a quasi-order.
\end{theorem}

Much like $\mathsf{S4}$-algebras are related to
Heyting algebras, $\mathsf{S4}_{t}$-algebras are 
related to special Heyting algebras, known as \textit{bi-Heyting algebras}:

\begin{definition}
    A \textit{bi-Heyting algebra} is a lattice $H$ such that both $(H,\leq)$ and $(H,\geq)$ are Heyting algebras. We denote the implication in $(H,\leq)$ by $\rightarrow$ and the implication in $(H,\geq)$ by $\dot{\text{---}}$.
\end{definition}

\begin{remark}
    The binary operation $\dot{\text{---}}$ is usually referred to as the {\em co-implication} and is sometimes denoted by $\leftarrow$ \cite{Wolter1998coimplication}. We follow \cite{Rauszer1974semiboolean} in denoting it by $\dot{\text{---}}$. 
\end{remark}

The original motivation for the study of such algebras concerns the ``problem of dualism" posed by McKinsey and Tarski \cite[Appendix]{mckinseytarskiclosedelements}.
They were studied in detail 
by Rauszer \cite{Rauszer1974semiboolean}, Esakia \cite{esakia1975problem}, Wolter \cite{Wolter1998coimplication}, and others. 
In particular, the Esakia representation for Heyting algebras takes on the following form for bi-Heyting algebras. Call an Esakia space $(X,\le)$ a \textit{bi-Esakia space} if ${\uparrow}U$ is clopen for each clopen $U \subseteq X$.

\begin{theorem} [Esakia representation for bi-Heyting algebras]
    Each bi-Heyting algebra is isomorphic to the bi-Heyting algebra of clopen upsets of a bi-Esakia space.
\end{theorem}

Note that in a bi-Esakia space $X$, given two clopen upsets $U$ and $V$, we have 
\begin{equation*}
    U \dot{\text{---}} V = {\uparrow}(V \, \setminus \, U). 
\end{equation*}

Like Esakia spaces are the partially ordered $\sf S4$-spaces, bi-Esakia spaces are the partially ordered $\mathsf{S4}_{t}$-spaces. The relationship between $\mathsf{S4}_{t}$-algebras and bi-Heyting algebras mirrors perfectly that of $\mathsf{S4}$-algebras and Heyting algebras:

\begin{samepage}
\begin{proposition}\
\begin{enumerate}
    \item If $(B,\Box_{F},\Box_{P})$ is an $\mathsf{S4}_{t}$-algebra then $B_{\Box_{F}}$ is a bi-Heyting algebra.
    \item If $(X,R)$ is an $\mathsf{S4}_{t}$-space, then $(\rho(X),\leq)$ is a bi-Esakia space.
\end{enumerate}
\end{proposition}
\end{samepage}

Note that if $X$ is bi-Esakia, then both $(X,\leq)$ and $(X,\geq)$ are Esakia spaces. Therefore, 
as an immediate consequence of Proposition \ref{Dual version of Johnstones theorem}, we obtain: 

\begin{proposition}[ZF+BPI]
    If $X$ is a bi-Esakia space, then $\max X\neq \emptyset$ and $\min X\neq \emptyset$.
\end{proposition}

As a result, all global maximality and minimality principles for bi-Heyting algebras are equivalent to BPI. This is reminiscent of what happens in 
Boolean algebras. 
However, the situation changes when we relativize to nonempty closed subsets:
\begin{samepage}
\begin{definition}\
\begin{enumerate}
    \item The \textit{bi-Esakia Maximality Principle}:
    \begin{align*}
    \tag{biREMax}
    &\text{ If $X$ is a bi-Esakia space and $C\subseteq X$ is a}\\
    &\text{nonempty closed subset, then $\max C\neq \emptyset$.}
\end{align*}
    \item \textit{Fine's Tense Maximality Principle}:
    \begin{align*}
    \tag{FtMax}
    &\text{ If $(X,R)$ is an $\mathsf{S4}_{t}$-space and $C\subseteq X$ is a }\\
    &\text{nonempty closed subset, then $\mathsf{qmax}C\neq \emptyset$.}
\end{align*}
\end{enumerate}
\end{definition}
\end{samepage}
 
\begin{theorem}[ZF+BPI]\label{Equivalence of Esakias bi-Maximality with the Axiom of Choice}
    The following principles are equivalent:
    \begin{enumerate}
        \item \textnormal{biREMax}; 
        \item \textnormal{FtMax}; 
        \item \textnormal{AC}. 
    \end{enumerate}
\end{theorem}
\begin{proof}
Since AC implies FMax, it also implies FtMax. Moreover, since every bi-Esakia space is an $\mathsf{S4}_{t}$-space, FtMax implies biREMax. It remains to show that biREMax implies AC, for which it is sufficient to show that biREMax implies PMax.
For this,
the key observation is that in the proof of Theorem \ref{Equivalence of Esakias Maximality with the Axiom of Choice}, 
$\mathbf{2}^{Z}$ is not only an Esakia space, but 
a bi-Esakia space, given that $\mathbf{2}$ is a bi-Esakia space and products of bi-Esakia spaces are bi-Esakia 
(since products are computed pointwise and for bi-Esakia spaces, both the space and its order-dual are Esakia spaces). Consequently, we can embed an arbitrary nonempty Priestley space $X$ in 
$\mathbf{2}^{Z}$, 
and since the latter is a bi-Esakia space, conclude that $\max X\neq \emptyset$. 
\end{proof}

\begin{samepage}
\begin{remark}\
\begin{enumerate}
        \item While we have focused here on $\mathsf{S4}_{t}$, one could obtain similar results for $\mathsf{K4}_{t}$ and even for weakly transitive tense logics, using the same ideas and techniques as in Section \ref{Fines Maximality Principle} (reflexivization) and Section \ref{Maximality in weaker modal systems} (reflexive and transitive closure).
        \item Similar to 
        Remark \ref{rem: WTMax algebraically}, we could consider weaker versions of 
        biREMax and FtMax by relativizing to clopen subsets. These are derivable from BPI. The corresponding algebraic formulations, 
        using relativizations, would be equivalent to BPI. 
\end{enumerate}
\end{remark}
\end{samepage}

The algebraic and topological principles considered in this paper are summarized below, together with their equivalences with AC or BPI, emphasizing the base theory in which these equivalences were established.

\vspace{3mm}

\begin{center}
    \begin{tabular}{|c|c|}
    \hline
    \multicolumn{2}{|c|}{ Algebraic Global Principles }\\
    \hline
          ScMax & Scott Maximality  \\
          \hline
         DLMax & Distributive Lattice Maximality for Ideals\\
         \hline
         HMax & Heyting Maximality for Ideals\\
         \hline
         LMax & Locale Maximality for Ideals\\
         \hline
         SLMax & Spatial Locale Maximality for Ideals\\
         \hline
         CLMax & Coherent Locale Maximality for Ideals\\
         \hline
    \end{tabular}
\end{center}

\begin{center}
    \begin{tabular}{|c|c|}
    \hline
    \multicolumn{2}{|c|}{ Algebraic Relativized Principles }\\
    \hline
          RHMax & Relativized Heyting Maximality for Filters  \\
          \hline
         FAMax & Fine's Algebraic Maximality\\
         \hline
         FAMax$^{\mathsf{S4}}$ & Fine's Algebraic Maximality for $\mathsf{S4}$\\
         \hline
    \end{tabular}
\end{center}

\begin{center}
    \begin{tabular}{|c|c|}
    \hline
    \multicolumn{2}{|c|}{ Topological Global Principles }\\
    \hline
          PMax & Priestley Maximality  \\
          \hline
         PMin & Priestley Minimality\\
         \hline
    \end{tabular}
\end{center}

\begin{center}
    \begin{tabular}{|c|c|}
    \hline
    \multicolumn{2}{|c|}{ Topological Relativized Principles }\\
    \hline
          REMax & Esakia Maximality  \\
          \hline
         FMax & Fine Maximality\\
         \hline
         FMax$^{\mathsf{S4}}$ & Fine's Maximality for $\mathsf{S4}$\\
         \hline
         FSMax & Fine's Strong Maximality\\
         \hline
         FSMax$^{S4}$ & Fine's Strong Maximality for $\mathsf{S4}$\\
         \hline
         WTMax & Weakly Transitive Maximality\\
         \hline
    FtMax & Fine's Tense Maximality\\
         \hline
         biREMax & bi-Esakia Maximality\\
         \hline
    \end{tabular}
\end{center}

\begin{figure}[h]
    \centering

    \begin{tikzpicture}
    \node at (4,1) {\textbf{Equivalents of AC
    ($\mathsf{ZF})$}};
        \node at (1,0) {ScMax};
        \node at (1,-1) {LMax};
        \node at (3,0) {DLMax};
        \node at (3,-1) {SLMax};
        \node at (5,0) {HMax};
        \node at (5,-1) {CLMax};
        \node at (7,0) {RHMax};
        \node at (3,-6) {FAMax};
        \node at (5,-6) {FAMax$^{\mathsf{S4}}$};
        \node at (4,-2) {\textbf{Equivalents of AC}
        ($\mathsf{ZF+BPI}$)};
        \node at (1,-3) {PMax};
        \node at (1,-4) {FSMax};
        \node at (3,-3) {PMin};
        \node at (3,-4) {FSMax$^{\mathsf{S4}}$};
        \node at (5,-3) {REMax};
        \node at (5,-4) {FtMax};
        \node at (7,-3) {WTMax};
        \node at (7,-4) {biREMax};
        \node at (4,-5) {\textbf{Equivalents of BPI ($\mathsf{ZF}$)}};
        
        \draw (0,0.5) -- (8,0.5) -- (8,-1.5) -- (0,-1.5) -- (0,0.5);

        \draw (0,-2.5) -- (8,-2.5) -- (8,-4.5) -- (0,-4.5) -- (0,-2.5);
         \draw (2,-5.5) -- (6,-5.5) -- (6,-6.5) -- (2,-6.5) -- (2,-5.5);
    \end{tikzpicture}    \caption{Equivalences of Maximality Principles}
    \label{fig:implicatinsbetweenprinciples}
\end{figure}


\printbibliography

\end{document}